\documentclass[12pt]{amsart}
\usepackage{amssymb,amscd,amsfonts,mathabx,mathrsfs}
\usepackage[arrow,matrix,curve]{xy}
\usepackage{fullpage}
\usepackage{graphicx}
\usepackage{slashed}
\renewcommand{\baselinestretch}{1.1}

\newcommand\dirac{\partialslash}

\newcommand\RE{\mathbb{R}}
\newcommand\R{\mathbb{R}}

\newcommand\ZA{\mathbb{Z}}
\newcommand\Z{\mathbb{Z}}

\newcommand\C{\mathbb{C}}

\usepackage{color}
\definecolor{darkgreen}{cmyk}{1,0,1,.2}
\definecolor{m}{rgb}{1,0.1,1}
\definecolor{green}{cmyk}{1,0,1,0}

\definecolor{test}{rgb}{1,0,0}
\definecolor{cmyk}{cmyk}{0,1,1,0}

\newcommand\Index{\operatorname{Index}}
\newcommand\Tr{\operatorname{Tr}}

\newcommand\Rank{\operatorname{Rank}}

\newcommand\tr{\operatorname{tr}}

\newcommand\E{\mathcal E}
\newcommand\maE{\mathcal E}

\newcommand\cS{\mathcal S}
\newcommand\maS{\mathcal S}

\theoremstyle{plain}
\newtheorem{theorem}{Theorem}[section]
\newtheorem{lemma}[theorem]{Lemma}
\newtheorem{proposition}[theorem]{Proposition}
\newtheorem{corollary}[theorem]{Corollary}

\theoremstyle{definition}
\newtheorem{definition}[theorem]{Definition}
\newtheorem{definition*}{Definition}

\theoremstyle{remark}
\newtheorem{remark}[theorem]{Remark}

\newtheorem{remarks*}{Remarks}

\begin{document}

\title[IPositive scalar curvature and twisted Dirac complexes]
{Conformal invariants of twisted Dirac operators\\ and positive scalar curvature}
\footnote{\today}
\author{Moulay Tahar Benameur}
\address{Laboratoire et D\'epartement de Math\'ematiques,
 UMR 7122,
 Universit\'e de Metz et CNRS,
 B\^at. A, Ile du Saulcy,
 F-57045 Metz Cedex 1, 
 France}
\email{benameur@math.univ-metz.fr}

\author{Varghese Mathai}
\address{Department of Mathematics, University of Adelaide,
Adelaide 5005, Australia}
\email{mathai.varghese@adelaide.edu.au}

\begin{abstract}
For a closed, spin, odd dimensional Riemannian manifold $(Y,g)$, 
we define the rho invariant $\rho_{spin}(Y,\E,H, [g])$ for the twisted Dirac operator $\dirac^\E_H$
on $Y$, acting on sections of a flat hermitian vector bundle $\E$ over $Y$,
where $H = \sum i^{j+1} H_{2j+1} $ is an odd-degree closed differential form on $Y$ 
and $H_{2j+1}$ is a real-valued differential form of degree ${2j+1}$. We prove that it only depends on
the conformal class $[g]$ of the metric $g$.
In the special case when $H$ is a closed 3-form, we use a Lichnerowicz-Weitzenb\"ock 
formula for the square of the twisted Dirac operator,  which in this case has 
no first order terms, to show that 
 $\rho_{spin}(Y,\E,H, [g])=\rho_{spin}(Y,\E, [g])$ for all $|H|$ small enough, whenever $g$ is conformally equivalent to a 
 Riemannian metric of positive scalar curvature.
When $H$ is a top-degree form on an oriented three dimensional manifold, 
we also compute  $\rho_{spin}(Y,\E,H)$. 
\\
\end{abstract}



\date{}
\maketitle

\section*{Introduction}

In an earlier paper \cite{BM}, we extended some of the results of
Atiyah, Patodi and Singer \cite{APS1, APS2, APS3} 
on the signature 
operator on an oriented, compact manifold with boundary, to the twisted case.
Atiyah, Patodi and Singer also studied the {\em Dirac operator} $\dirac^\E$ on an odd dimensional, closed, spin manifold,
which is self-adjoint and elliptic, having spectrum in the real numbers. For this 
(and other elliptic self-adjoint operators), they introduced the eta invariant
which measures the spectral asymmetry of the operator and is a 
spectral invariant. Coupling with flat bundles, they introduced the closely related rho invariant,
which has the striking property that it is independent of the choice of Riemannian
metric needed in its definition, when it is reduced modulo $\ZA$. In this paper we generalize the construction of 
Atiyah-Patodi-Singer to the twisted Dirac operator  $\dirac^\E_H$ with a closed, odd-degree differential 
form as flux and with coefficients in a flat vector bundle.

More precisely, let $X$ be a $2m$-dimensional compact, spin Riemannian manifold without boundary, $\E$
a flat hermitian vector bundle over $X$ and $H$ a closed, odd degree differential form on $Y$. Consider the twisted Dirac operator $\dirac^\E_H=c\circ \nabla^{\E, {H}} = \dirac^\E + c(H)$ 
where $c$ denotes Clifford multiplication, $\nabla^{\E, H}$ denotes the flat superconnection 
$\nabla_X^\E+c({H})$, and
$\nabla^\E$ is the canonical flat hermitian connection on $\E$.
Then $\dirac^\E_H$ anticommutes with the usual grading involution and it is  self-adjoint if and only if 
\begin{equation}\label{eqn:fluxform}
H = \sum i^{j+1} H_{2j+1} 
\end{equation}
where $H_{2j+1}$ are  real-valued differential form 
of degree ${2j+1}$. It is only in this case that one gets a { {generalization of the usual}} Dirac 
operator on $X$, in contrast to the 
case of the twisted de Rham complex, cf . \cite{BCMMS,MS03,AS,MW,MW2,RW}. 

{ {When the compact spin manifold $X$ has non empty boundary and }}assuming that the 
Riemannian metric is of product type near the boundary { {and that $H$ satisfies the absolute boundary condition}}, we explicitly identify 
the twisted Dirac operator near the boundary
to be 
$\dirac^{\tilde\E}_{\tilde H} = \sigma\left(\frac{\partial}{\partial r} + \dirac^\E_H \right),$ where 
$r$ is the coordinate in the normal direction, $\sigma$ is a bundle isomorphism and 
finally, the self-adjoint elliptic operator  given on $\Omega^{2h}(\partial X, \E)$ by
$\dirac^\E_H$ is the Dirac operator on the boundary. 
$\dirac^\E_H$ is an elliptic self-adjoint operator, and 
by \cite{APS1}, the non-local boundary condition given by $P^+(s\big|_{\partial X})=0$
where $P^+$ denotes the orthogonal projection onto the eigenspaces with positive eigenvalues,
makes the pair $(\dirac^\E_H; P^+)$
into an elliptic boundary value problem. Applying the Atiyah-Patodi-Singer index theorem, and computing the local contribution when $H$ is closed, we get 
$$
\Index(\dirac^{\tilde\E}_{\tilde H}; P^+) =  \Rank (\E) \int_Y \alpha_0(H) -\frac{1}{2} \left(\dim(\ker(\dirac^\E_H))
+\eta(\dirac^\E_H)\right),
$$
where $\eta(\dirac^\E_H)$ denotes the eta invariant of the self-adjoint operator $\dirac^\E_H$ and $\alpha_0(H)$ is a differential form which does not depend on the flat hermitian bundle $\E$.
This is the main tool used to prove our results about the twisted Dirac rho invariant, defined below.

Let now $Y$ be a closed, oriented, { {$(2m - 1)$}}-dimensional Riemannian spin manifold and 
$H = \sum i^{j+1} H_{2j+1} $ an 
odd-degree, closed differential form on $Y$ where $H_{2j+1}$ is a real-valued differential form 
of degree ${2j+1}$. Denote by $\E$ a hermitian flat vector bundle over $Y$ with the 
canonical flat connection $\nabla^\E$.
Consider the twisted Dirac operator $\dirac^\E_H=c\circ \nabla^{\E,  H}= \dirac^\E + c(H)$
Then $\dirac^{ {\E}}_H$ is a self-adjoint
elliptic operator and let $\eta(\dirac^\E_H )$ denote its eta invariant. The twisted Dirac rho invariant 
$\rho_{spin}(Y, \E, H, [g])$ is defined to be the difference 
$$
\rho_{spin}(Y, \E, H, [g]) = \frac{1}{2} \left(\dim(\ker(\dirac^\E_H))
+\eta(\dirac^\E_H)\right)
 - \Rank (\E) \frac{1}{2} \left(\dim(\ker(\dirac_H))
+\eta(\dirac_H)\right),
$$
where 
$\dirac_H$ is the same twisted Dirac operator corresponding to the trivial line bundle.
Although the eta invariant  $\eta(\dirac^\E_H )$  is a priori only a spectral invariant, we show that the  twisted Dirac rho invariant, $\rho_{spin}(Y,\E,H, [g])$ depends only on the conformal class of the Riemannian metric $[g]$.
We compute it for 3-dimensional manifolds with a 
degree three flux form, Corollary \ref{cor:3dcalc}. This is done via the important method of spectral flow as in \cite{APS3,Getzler}. In section \ref{sec:pos}, we analyse the special case when $H$ is a closed 3-form, using 
a Lichnerowicz-Weitzenb\"ock 
formula for the square of the twisted Dirac operator,  which in this case has 
no first order terms, to show that 
%whenever $X$ is a closed spin manifold of positive scalar curvature, then
 $\rho_{spin}(Y,\E,H, [g])=\rho_{spin}(Y,\E, [g])$ for all $|H|$ small enough, whenever $g$ is conformally equivalent to a 
 Riemannian metric of positive scalar curvature.

The twisted analogue of analytic torsion was studied in \cite{MW, MW2,MW3} and 
is another source of inspiration for this paper. We mention that twisted Dirac operators,
known as {\em cubic Dirac operators} have been studied in representation theory
of Lie groups on homogeneous spaces \cite{Kostant, Slebarski}. It  also appears in the study of Dirac operators on 
loop groups and their representation theory \cite{Landweber}.

\medskip

\noindent {\bf Acknowledgments}
M.B. ~acknowledges support from the CNRS institute INSMI especially from the PICS-CNRS ``Progress in Geometric Analysis and Applications''.
V.M.~acknowledges support from the Australian Research Council.
Part of this work was done while the second author was visiting the ``Laboratoire de Math\'ematiques et Applications'' at Metz.

\section{Twisted Dirac operator and Lichnerowicz-Weitzenb\"ock formulae}

\subsection{The twisted analog of the Dirac operator}

Assume that $X$ is an oriented manifold of dimension $n=2m$ and that $H$ is  given complex differential form of positive degree on $X$. We denote by $\check{H}$ and ${\hat H}$ the differential forms
$$
\check{H}=\sum_{k\geq 1} H_k/k \text{ and } {\hat H}= \sum_{k\geq 1} k H_k, \quad \text{ if } H=\sum_{k\geq 1} H_k.
$$
We assume  that $X$ is endowed with a spin structure and denote by $\maS=\maS^+\oplus \maS^-$ the $\Z_2$-graded spin bundle. We fix a  unitary-flat hermitian bundle $\maE, \nabla^\E$ and denote by  $\gamma$ the grading involution obtained on $\maS\otimes \E$.. 
{ {
We are interested in the twisted Dirac operator $\dirac^\E_H = c\circ \nabla^{\E, \check{H}}$ 
where 
$$
\nabla_X^{\E, \check{H}}=\nabla_X^\E+i_X \check{H}\,\cdot\,,
$$ 
and $\nabla^\E$ is the canonical flat hermitian connection on $\E$. So 
$$
\dirac^\E_H =  \dirac^\E + c(H),
$$
where $c(H)$ is the action of the differential form $H$ by Clifford multiplication on the Clifford module $\maS\otimes \E$.
Then we have,}}

\begin{proposition}
In the notation above, 
\begin{itemize}
\item $
\dirac^\E_H \gamma = - \gamma \dirac^\E_H \Longleftrightarrow H\in \Omega^{odd} (X, \C).
$
\item The twisted Dirac operator $\dirac^\E_H$ is self-adjoint if and only if 
$$
H = \sum_{k\geq 1} (A_{4k} + i B_{4k-2}) + (A_{4k-1} + i B_{4k-3}),
$$ 
where the differential forms $A_j$ and $B_j$ are real differential forms of degree $j$ on $X$.
\end{itemize}
\end{proposition}

So, in the even dimensional case, it is  only in the  case $H\in \Omega^{4\bullet +3} + i \Omega^{4\bullet + 1}$ that one gets { {a twisted version of the usual}} Dirac operator. When the ambient manifold is odd dimensional, then we may consider also even forms. Compare with \cite{BCMMS,AS,MW,MW2}. Notice that this condition coincides with (and explains) the one obtained for the twisted signature operator in \cite{BM}. 

\begin{proof} We have $\dirac^\E_H = \dirac^\E + c(H)$, where  $\dirac^\E $ is the usual Dirac operator acting on sections of
$\cS(X)\otimes \E$. It is a classical result that  $\dirac^\E \gamma = - \gamma \dirac^\E$ and that $\dirac^\E$ is self-adjoint. Recall that Clifford multiplication by a differential form of degree $j$ is self-adjoint if and only if $j$ is congruent to $0$ or $3$ modulo $4$ and is skew adjoint otherwise. This proves the second item.

Recall now that $\gamma$ is given locally using an orthonormal basis by Clifford multiplication with the differential volume form
$$
\omega = i^m e_1\cdots e_{2m}.
$$
A straightforward local computation then shows that for any complex differential form $\alpha$ of degree $j$, we have $
c(\alpha)\circ \gamma - (-1)^j \gamma\circ c(\alpha) = 0.$
\end{proof}

\subsection{Lichnerowicz-Weitzenb\"ock formulae for odd degree twist}

Let $H \in \Omega^{odd}(X)$ be a closed differential form of odd degree 
and $c(H)$ be the image of $H$ in the sections of the Clifford algebra bundle ${\rm Cliff}(TX)$. 
Then $\nabla^\E + H$ is a superconnection on the trivially graded, flat bundle $\E$ over $X$. Then 
$ (\nabla^\E + H)^2 =dH = 0$
is the curvature of the superconnection which is flat. 

Let $\dirac^\E$ denote the Dirac operator acting on $\E$-valued spinors on $X$. If $\{e_1, \ldots, e_n\}$ is a local orthonormal 
basis of $TX$, then we have the expression,
$$
\dirac^\E = \sum_{j=1}^n c(e_j) \nabla^\E_{e_j}
$$
Let $R^g$ denotes the scalar curvature of the Riemannian metric. Then as shown in \cite{Bismut} the spinor Laplacian
$$\Delta_H^\E = - \sum_{j=1}^n \left(\nabla^\E_{e_j} + c(\iota_{e_j} H)\right)^2$$ is a positive operator that 
does not depend on the local orthonormal basis  $\{e_1, \ldots, e_n\}$
of $TX$. Here $\iota_{e_j}$ denotes contraction by the vector $e_j$.

Then the following is a consequence of Theorem 1.1 in \cite{Bismut}.

\begin{theorem}[Lichnerowicz-Weitzenb\"ock formulae\cite{Bismut}]\label{thm:bismut}
Let $H$ be a closed, odd degree differential form on $Y$.
Then the following identities hold,
$$
\left(\dirac^\E_{H}\right)^2 = \Delta_H^\E  + \frac{R^g}{4} 
+ c(H)^2 +  \sum_{j=1}^n c(\iota_{e_j} H)^2,
$$
where $R^g$ denotes the scalar curvature of the Riemannian spin manifold $X$
and the last 2 terms on the right hand side satisfy,
$$
 c(H)^2 +  \sum_{j=1}^n c(\iota_{e_j} H)^2  = \sum_{j_1<j_2\ldots<j_k, \, k\ge 2}
 (-1)^{\frac{k(k+1)}{2}} (1-k) c((\iota_{e_{j_1}}\iota_{e_{j_2}}\ldots \iota_{e_{j_k}}H)^2).
$$
\end{theorem}

As a corollary of Theorem \ref{thm:bismut}, one has the following 
%well known result (Theorem 1.3 \cite{Bismut} and Theorem 6.2 of \cite{AF}), one has the following 
special Lichnerowicz-Weitzenb\"ock formula

\begin{theorem}\label{thm:weit}\cite{Bismut,AF}
Let $H$ be a closed differential 3-form. 
Then 
$$
\left(\dirac^\E_{H}\right)^2 = \Delta_H^\E  +  \frac{R^g}{4} - 2 \vert H\vert^2
$$
where $R^g$ denotes the scalar curvature of the Riemannian spin manifold $X$ and $|H|$ the length of $H$.
\end{theorem}

\section{Eta invariants of twisted Dirac operators}
\label{sect:etarho}

In this section, we define rho invariant $\rho_{spin}(X,\E,H)$ and eta invariant $\eta(\dirac^\E_H)$
of the twisted Dirac operators. But let us review the boundary Dirac operator first.

Let $X$ be an $2m$ dimensional compact  spin manifold with boundary $Y$, $\E$
a flat hermitian vector bundle over $X$ and $H$ is a closed, odd degree differential form on $X$ such that 
$$
H\in \Omega^{4\bullet+3} (X, \R) \oplus i \Omega^{4\bullet +1} (X, \R).
$$
 Denote by $\maS$ the spin bundle and 
consider as above the involution 
$\gamma_X=i^m c_X(e_1\cdots e_{2m})\otimes 1_\E$ on $\maS \otimes \E$-sections over $X$.  According to this involution we decompose as usual $\maS \otimes \E$ into 
$$
\maS \otimes \E = \maS^+\otimes  \E \oplus \maS^-\otimes \E.
$$
We may and shall assume that the local orthonormal basis $(e_k)_{1\leq k \leq 2m}$ 
is decomposed on the boundary into the local orthonormal basis $(e_k)_{1\leq k\leq 2m-1}$ and the inward unit vector $e_{2m}=\partial/\partial r$ which is orthogonal to 
the boundary. Clifford multiplication by $e_{2m}$ is denoted below by $\sigma$. We shall then also consider the self-adjoint involution $\gamma_Y= i^{m} c(e_1\cdots e_{2m-1})$.

{ { We denote by $i^*H$ the restriction of $H$  to the boundary $Y=\partial X$. We shall say that $H$ is boundary-compatible if there exists a collar neighborhood 
$p_\epsilon: X_\epsilon\cong (-\epsilon, 0] \times Y \to Y$ such that }}
\begin{equation}\label{bdry.compatible}
{ {H  \vert_{X_\epsilon}  = p_\epsilon^* (i^*H).}}
\end{equation}

As usual the spinor bundle $\maS_Y$ on the boundary manifold $Y$ for the induced spin structure from the fixed one on $X$ is naturally identified with $\maS^+$ so that the Clifford representations are related by:
$$
c_Y (U) = -\sigma \cdot c_X(U), \quad \text{ for } U\in TY\subset TX.
$$

\begin{lemma}\label{lem:bdry-iden}
Under the above assumptions, suppose further that the 
riemannian metric on $Y$ and the hermitian metric on $\E$ are of product type near the boundary { {and that the closed, odd degree differential form $H$ is  boundary-compatible}}.
We identify as usual the restriction of $\maS^+$ to the boundary $Y$ with the spin bundle $\maS_Y$. Then, in the collar neighborhood $X_\epsilon$, we explicitly identify 
the twisted Dirac operator $\dirac_H^{X,\E}$ with the operator
$$
\sigma \left(\frac{\partial}{\partial r} + \dirac^{Y, \E\vert_Y}_{H\vert_Y} \right),
$$ 
where 
 $\dirac^{Y, \E\vert_Y}_{H\vert_Y}$ is the self-adjoint elliptic twisted Dirac operator on $Y$ defined as before by
$
\dirac^{Y, \E\vert_Y}_{H\vert_Y} =  c\circ \nabla^{\E, \check{H}_{\partial Y}}
$
acting on $\cS(\partial Y)\otimes \E$.
\end{lemma}

\begin{proof}

Choose an orthonormal basis $(e_1, \cdots, e_{2m})$ as above near a point on the boundary $Y$ such that $e_{2m}=\frac{\partial}{\partial r}$, 
$\sigma$ anti commutes with the grading involution $\gamma_X$ and is itself a well defined involution of $\maS\vert_{X_\epsilon}$. Given a spinor 
$\phi\in \Gamma (Y, \maS_Y\otimes \E_Y)$ over the boundary manifold $Y$, the pull-back section $\pi_\epsilon^* \phi$ is identified with a section 
of $\maS^+\otimes \E \cong p_\epsilon ^*(\maS_Y\otimes \E_Y)$. Therefore, for any smooth function $f$ on $(-\epsilon, 0]$, we compute
\begin{eqnarray*}
 \sigma\cdot \dirac^{X, \E}_H (f\pi_\epsilon^* \phi) & = & f  \sum_{i=1}^{2m-1} \sigma c_X(e_i) \pi_\epsilon^*(\nabla^Y_{e_i} \phi) + \frac{\partial f}{\partial r}
\sigma^2 \pi_\epsilon^*\phi + f \sigma (H\cdot \pi_\epsilon^*\phi)\\
& = &  -  f \sum_{i=1}^{2m-1} c_Y(e_i) \pi_\epsilon^*(\nabla^Y_{e_i} \phi) - \frac{\partial f}{\partial r} \pi_\epsilon^*\phi + f \sigma c_X(H)(\pi_\epsilon^*\phi)\\
& = &  -  f \pi_\epsilon^* \left(\sum_{i=1}^{2m-1} c_Y(e_i) \nabla^Y_{e_i} \phi\right) - \frac{\partial f}{\partial r} \pi_\epsilon^*\phi - f \pi_\epsilon^* (c_Y(H)(\phi))\\
& = & - f \pi_\epsilon^* (\dirac^{Y, \E\vert_Y, H\vert_Y}\phi) - \frac{\partial f}{\partial r} \pi_\epsilon^*\phi. 
\end{eqnarray*}
Composing again by $\sigma$ and using $\sigma^2 = -1$, we conclude.
\end{proof}

Let us briefly recall the definition of the eta invariant. 
Given a self-adjoint elliptic differential operator $A$ of order $d$ on a closed 
oriented manifold $Y$ of dimension $2m-1$,
the {\em eta-function} of $A$ was defined in \cite{APS2} as
$$
\eta(s,A):=\Tr'(A|A|^{-s-1}), 
$$
where $\Tr'$ stands for the trace restricted to the subspace orthogonal to
$\ker(A)$. By \cite{APS1, APS2, APS3}, $\eta(s,A)$ 
is holomorphic when $\Re(s)>2m-1/d$ and can
be extended meromorphically to the entire complex plane with possible simple
poles only.   
It is related
to the heat kernel by a Mellin transform
$$
\eta(s,A)=\frac{1}{\Gamma(\frac{s+1}{2})}\int_0^\infty t^{\frac{s-1}{2}}\Tr(Ae^{-t{\widetilde A}^2})\,dt.
$$

The eta function is then
known to be holomorphic  at $s=0$, more precisely   by  \cite{APS3} one has:

\begin{theorem}[\cite{APS3}]\label{zeta-hol}
Let $A$ be a first order self-adjoint elliptic operator on the odd dimensional closed manifold $X$. Then the eta function $\eta(s,A)$ of $A$ has a meromorphic continuation 
to the complex plane with no pole at $s=0$. The eta invariant $\eta(A)$ of $A$ is then defined as $\eta(0,A)$.
\end{theorem}

The {\em eta-invariant}
of $A$ is thus defined as
$$
\eta(A) = \eta(0, A).
$$

Consider now our twisted Dirac operator 
$\dirac^\E_H $ defined in the previous subsection for general odd dimensional spin manifold $Y$ and flat bundle $\E\to Y$, and which is a self-adjoint
elliptic differential operator and let $\eta(\dirac^{Y, \E}_H )$ then denote its well defined eta invariant. 

Following \cite{APS2}, we set 
\begin{equation}\label{eqn:xi}
\xi_{spin} (Y, \E, H):= \frac{1}{2} \left( \dim(\ker(\dirac^\E_H)) + \eta(\dirac^\E_H) \right).
\end{equation}

\begin{definition}
The {\em twisted Dirac rho invariant} $\rho_{spin}(Y, \E, H)$ is defined to be 
$$
\rho_{spin}(Y, \E, H) := \xi_{spin} (Y, \E, H ) - {\rm rank}(\E)\, \xi_{spin} (Y, H ),
$$
 where $\xi_{spin} (Y, H )$ is the invariant $\xi$ corresponding to the case where the flat hermitian bundle $\E$ is the trivial line bundle.
\end{definition}

As for the untwisted case, we denote this reduction of $\rho_{spin}(Y, \E, H)$ mod $\ZA$ by $
\bar\rho_{spin} (Y, \E, H)$. Then, 
the reduced twisted  rho invariant $\bar\rho_{spin}(Y,\E,H)$ is independent of the
choice of the Riemannian metric on $X$ and the Hermitian metric on $\E$. It is also a cobordism invariant of the triple $(Y,\E,H)$. In the case of positive scalar curvature, we shall though be able to work  with the real invariant $\rho_{spin} (Y, \E, H)$. 
These results will be established in section \ref{sect:var}.

\subsection{The twisted Dirac index for manifolds with boundary}
\label{sect:twistedsig}

The goal of this section is to review  the Atiyah-Patodi-Singer index theorem for the 
twisted Dirac operator $\dirac^{X, \E}_{H}$ 
with non-local boundary conditions. Here and as before, $\E$ is a unitary flat hermitian bundle on $X$ and $H$ is a closed,  odd degree differential form on $X$ which is boundary compatible
and belongs to $\Omega^{4\bullet + 3} (X, \R) \oplus i \Omega^{4\bullet + 1} (X, \R)$. The notation $H\vert_Y$ stands for the restriction $i^*H$ of $H$ to the boundary manifold 
$\partial X=Y$.
Proposition \ref{prop:sig1} below, can be deduced from \cite{APS1}.
%is one of the key tools that will be used to prove our results on the eta and rho invariants for the twisted odd Dirac
%operator in the following section.

$\dirac^{\tilde \E}_{H}$ is an elliptic self-adjoint operator, and 
by \cite{APS1}, the non-local boundary condition given by $P^+(s\big|_{\partial X})=0$
where $P^+$ denotes the orthogonal projection onto the eigenspaces with positive eigenvalues,
yields an elliptic boundary value problem $(\dirac^{\tilde \E}_{ H}; P^+)$. Recall that $\dim (X)= 2m$.
By the Atiyah-Patodi-Singer index theorem \cite{APS1, APS3} and its extension in \cite{Gilkey}, we have

\begin{proposition}\label{prop:sig1}
In the notation above,
$$
\Index(\dirac^{X, \E}_{ H}; P^+) =  \Rank(\E) \int_X \alpha_0( H)  - \xi_{spin} (Y,\E\vert_Y,  {H\vert_Y}).
$$
where $ \xi_{spin} (Y,\E\vert_Y,  {H\vert_Y})$ is 
%denotes the  $xi invariant of the Dirac operator $\dirac^{Y, \E\vert_Y}_{H\vert_Y}$ acting on the spinors of 
%the boundary spin manifold $Y$ 
as defined in equation  \eqref{eqn:xi} and $\alpha_0(H)$ is the local contribution given by the APS theorem.
\end{proposition}

\begin{remark} 
The precise form of $\alpha_0(x)$ is unknown for general $H$. 
However in the case when $H=0$, the local index theorem cf. \cite{BGV} establishes that the Atiyah-Singer
$\widehat A$-polynomial applied to the curvature of the Levi-Civita connection, wedged by the Chern character of the flat bundle 
$\E$, is equal to
 $\alpha_0(x)$ is equal to  times the rank of $\E$.
   In the case when degree of $H$ is equal to 3, it  follows from a result of Bismut \cite{Bismut} that the Atiyah-Singer
$\widehat A$-polynomial applied to the curvature of a Riemannian connection defined in terms of the 
Levi-Civita connection together with a torsion tensor determined by $H$, wedged by the Chern character of the flat bundle 
$\E$,  is equal to
 $\alpha_0(x)$ times the rank of $\E$. The proof that in general, one also gets
  $\alpha_0(x)$ times the rank of $\E$ is exactly as argued in the Appendix to \cite{BM}.
\end{remark}

\subsection{Conformal invariance of the twisted spin eta and rho invariant}
\label{sect:var}
Here we prove the twisted spin rho invariant $\rho_{spin}(X,\E,H, [g])$ depends only on the
conformal class $[g]$ of the Riemannian metric $g$ on $X$ and the Hermitian metric on $\E$ needed
in its definition. 
The proof relies on the index theorem for twisted Dirac 
operator for spin manifolds with boundary, established in Proposition \ref{prop:sig1}.
We also state and prove the basic functorial properties of the twisted spin rho invariant.\\

\noindent{\em Conformal variation of the Riemannian metrics.}
We assume that $X$ is a compact spin manifold of dimension  $(2m-1)$.
Let $g$ be a Riemannian metric on $X$ and $g_\E$, an Hermitian metric 
on $\E$.
%Let $Q_{\bar k}$ ($k=0,1$) be the orthogonal projection from (the completion
%of) $C^{\bar k}$ to $\ker(\Delta_{\bar k})$. 
Suppose that $g$ is deformed smoothly and conformally along a one-parameter
family  $t\in\RE$.
%, then the operators $\ast$, $\sharp$ and 
%$\Gamma=\ast\sharp=\sharp\ast$ (see \$\ref{sect:scalar}) all depend smoothly
%on $u$.
The Dirac operator $\dirac^{X, \E}_g$ acting on the bundle of spinors over 
$(X, g)$,  is conformally covariant according to \cite{LM}, where ${\hat g} = e^f g$, for $f$ a 
smooth function on $X$,
\begin{equation}\label{eqn:conformal}
\dirac^{X, \E}_{\hat g} = e^{m f}\circ \dirac^{X, \E}_{g}\circ  e^{-(m-1) f}.
\end{equation}
It follows that the twisted Dirac operator, $\dirac^{X, \E}_H$ is also conformally covariant 
with the same weights.
Then we have,

\begin{theorem}[Conformal invariance of the spin rho invariant]\label{thm:indept}
Let $Y$ be a compact, spin manifold of dimension $2m-1$, $\E$, a flat hermitian vector
bundle over $Y$, and $H = \sum i^{j+1} H_{2j+1} $  an 
odd-degree {{closed}}  differential forms on $Y$ and $H_{2j+1}$ is a real-valued differential form 
homogeneous of degree ${2j+1}$. 
Then the spin rho invariant $\rho_{spin}(Y,\E,H, [g])$ of the twisted Dirac operator
depends only on the conformal class of the Riemannian metric on $Y$.
% and on the conformal class
%of the Hermitian metric on $\E$. 
%It is a diffeomorphism invariant of $X, \E$ and $H$.
\end{theorem} 

\begin{proof}
Consider the manifold with boundary $X = Y \times [0, 1]$, where the boundary $\partial X
= Y \times \{0\} - Y \times \{1\}$. 
Choose a smooth function $a(t), \, t\in [0,1]$ such that $a(t)\equiv 0$ near $t=0$ and 
$a(t)\equiv 1$ near $t=1$. Consider the metric $h=e^{2a(t) f}(g + dt^2)$ on $X$, 
which is also of product type near the boundary,
and let $\pi : X \to Y$
denote projection onto the first factor.

\begin{figure}[h]
\includegraphics[height=2in]{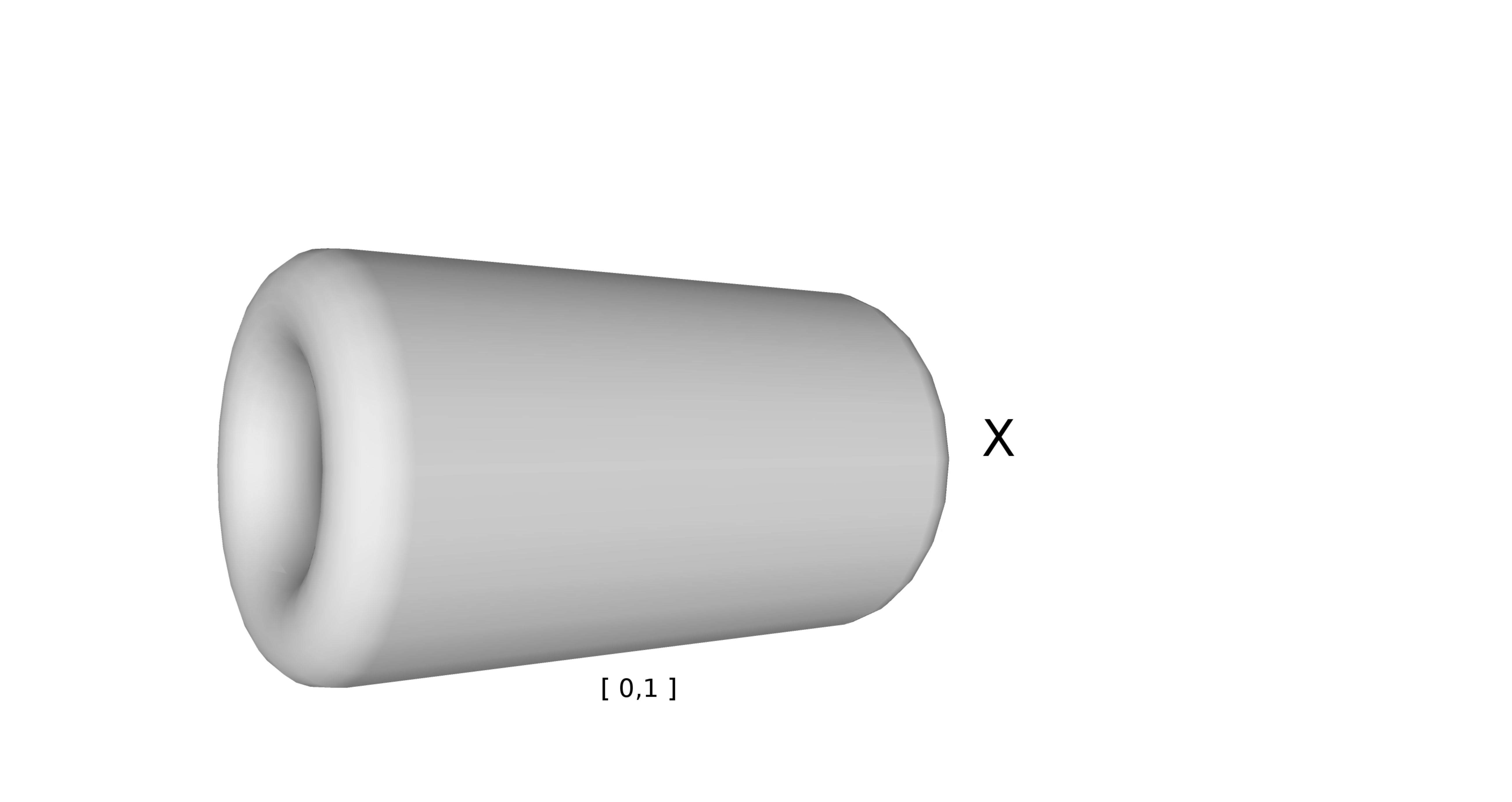}\\
\caption{}
\label{fig:box-metric}
\end{figure}

Applying the index theorem for the twisted Dirac operator, Proposition \ref{prop:sig1}, we get
\begin{equation}\label{eqn:sig1}
{\rm Index}\left(\dirac^{X, \pi^*(\E)}_{\pi^*(H)}, P^+_\E\right) =
\Rank(\E) \int_{X} \alpha_0^H +
\xi_{spin}(\dirac^\E_H, g)
-\xi_{spin}(\dirac^\E_H, \hat g).
\end{equation}
On the other hand, applying the same theorem to the trivial bundle ${\underline{\Rank (\E)}}$ of rank 
equal to $\Rank (\E)$, we get
\begin{equation}\label{eqn:sig3}
{\rm Index}
\left(\dirac^{X, {\underline{\Rank (\E)}}}_{\pi^*(H)}, P^+\right) 
=  \Rank(\E)\int_{X} \alpha_0^H + 
\Rank (\E) \left[\xi_{spin}(\dirac_H,  g)
-\xi_{spin}  (\dirac_H, \hat g)\right].
\end{equation}
Subtracting the equalities in \eqref{eqn:sig1} and \eqref{eqn:sig3} above, we get
\begin{equation}\label{eqn:sig5}
\rho_{spin}(Y, \E,H, g) - \rho_{spin}(Y,\E,H, \hat g) =
{\rm Index}\left(\dirac^{X, \pi^*(\E)}_{\pi^*(H)}, P^+_\E\right) - {\rm Index}\left(\dirac^{X, {\underline{\Rank (\E)}}}_{\pi^*(H)}, P^+\right).
\end{equation}
Each of the index terms on the right hand side of 
\eqref{eqn:sig5} has been shown in \cite{APS1} to be equal to the $L^2$-index of $\hat X$, which is 
$X$ together with infinitely long metric cylinders glued onto it at $\partial X$, plus a correction term (the dimension of the space of, limiting values of right handed spinors). The  $L^2$-index is a conformal invariant by \eqref{eqn:conformal}, and similarly, the correction term is also a conformal invariant. It follows by \eqref{eqn:sig5} that 
$$
\rho_{spin}(Y, \E,H, g) = \rho_{spin}(Y,\E,H, \hat g).
$$
\end{proof}

\section{Spectral flow and calculations of the twisted Dirac eta invariant}
\label{sect:calc}

\subsection{Spectral flow and the twisted Dirac eta invariant}
We employ the method of 
spectral flow for a path of self-adjoint elliptic operators, which is generically the net number of eigenvalues
that cross zero, which was first defined by \cite{APS3}, to get another formula for the difference 
$\eta(\dirac^\E_H) - \eta(\dirac^\E)$. Extensions of spectral flow to families using entire cyclic cohomology recently appeared in  \cite{BC}.
% which holds more generally than in Proposition \ref{prop:3mfd} above. 

Let $(Y, g)$ be an odd-dimensional  Riemannian spin manifold. 
Define $S^{-H} \in \Omega^1(Y, {\mathcal S}{\mathcal O}(TY))$, which is a degree one form with values in 
the skew-symmetric endomorphisms of $TY$, as follows. For $\alpha, \beta, \gamma \in TY$, set 
$$
g(S^{-H} (\alpha)\beta, \gamma) = - 2H(\alpha, \beta, \gamma).
$$
Let $\nabla^L$ denote the Levi-Civita connection of $(Y, g)$ and $\nabla^{-H} = \nabla^L + S^{-H} $ 
be the Riemannian connection whose curvature is denoted by $\Omega_Y^{-H}$.

\begin{proposition}\label{prop:getzler}
Let $Y$ be as before a $(2\ell-1)$-dimensional closed, spin, Riemannian manifold and
$\E$ is a flat vector bundle associated to an orthogonal or unitary
representation of $\pi_1(Y)$. Let $H\in \Omega^{3} (Y, \R)$ be a closed 
3-form such that
 the Dirac operators $\dirac^\E$ and $\dirac^\E_H$ are both invertible, then
$$
{\rm sf}(\dirac^\E, \dirac^\E_H) = \frac{\Rank(\E)}{(-2\pi i)^{\ell+1}} \int_Y H \wedge \widehat A(\Omega_Y^{-H}) 
 + \frac{1}{2}(\eta(\dirac^\E_H) - \eta(\dirac^\E)),
$$
where $\widehat A(\cdot)$ denotes the A-hat genus polynomial, $\Omega_Y^{-H}$ denotes the curvature of the 
Riemannian connection $\nabla^{-H}$ and $\, {\rm sf}(\dirac^\E, \dirac^\E_H)$ denotes the spectral flow of the smooth path of self-adjoint elliptic operators $\{\dirac^\E_{uH}, \, : \, u\in [0,1]\}$.
% is given by:
%$$
%2 \times {\rm sf}(\dirac^\E, \dirac^\E_H) = \eta(\dirac^\E_H) - \eta(\dirac^\E) 
%$$
\end{proposition}

\begin{proof} By theorem 2.6 \cite{Getzler}, 
$$
{\rm sf}(\widetilde \dirac^\E, \widetilde \dirac^\E_H) = \left(\frac{\epsilon}{\pi}\right)^{1/2}\int_0^1 \Tr\left(c(H)\circ
e^{-\epsilon {\dirac^\E_{uH}}^2}\right) du
+ \frac{1}{2}(\eta(\dirac^\E_H) - \eta(\dirac^\E))
$$
where $\widetilde \dirac^\E_H$ denotes the operator $ \dirac^\E_H + P_{\E, H}$ where $P_{\E, H}$ denotes the 
orthogonal projection onto the nullspace of $ \dirac^\E_H$. Following the idea of the proof of Theorem 2.8 in \cite{Getzler}
the technique of \cite{BGV}, and the Local Index Theorem 1.7 of Bismut \cite{Bismut} (here the fact that $H$ is a closed 3-form on $Y$ is used), we may make the following replacements:
\begin{align*}
\epsilon^{1/2}c(H)\qquad & \text{by} \qquad H\wedge\\  
e^{- (\epsilon^{1/2}\dirac^\E_{uH})^2} \qquad & \text{by} \qquad  \widehat A(\Omega_Y^{-H}) \wedge e^{(d+uH)^2}\\
\Tr(\cdot) \qquad & \text{by} \qquad \frac{-i\pi^{1/2}}{(-2\pi i)^{\ell+1}} \int_Y \tr(\cdot)
\end{align*}  
This then enables us to replace $\displaystyle \left(\frac{\epsilon}{\pi}\right)^{1/2}\int_0^1 \Tr\left(c(H)\circ
e^{-\epsilon {\dirac^\E_{uH}}^2}\right) du$ by 
$$
\frac{\Rank(\E)}{(-2\pi i)^{\ell+1}} \int_Y H \wedge \widehat A(\Omega_Y^{-H}),
$$
and the proposition follows.
\end{proof}

As a special case of  Proposition \ref{prop:getzler}, one obtains the following calculation,

\begin{corollary}\label{cor:3dcalc}
Let $Y$ is a compact spin Riemannian manifold of dimension 3 and
$\E$ is a flat vector bundle associated to an orthogonal or unitary
representation of $\pi_1(Y)$.
Let $H$ be a closed $3$-form on $Y$.
Consider the smooth path of self-adjoint elliptic operators $\{\dirac^\E_{uH}, \, : \, u\in [0,1]\}$
and assume that the Dirac operators $\dirac^\E$ and $\dirac^\E_H$ are both invertible. Then
$$
\eta(\dirac^\E_H) - \eta(\dirac^\E) = 2\, {\rm sf}(\dirac^\E, \dirac^\E_H)  + \frac{h}{2\pi^2}
$$
where $h=[H] \in \mathbb R\cong H^3(Y, \RE)$.
\end{corollary}

\section{Positive scalar curvature}\label{sec:pos}

Here we give applications of our results to closed spin manifolds that admit a Riemannian  metric of 
positive scalar curvature. 
Foundational works on metrics of positive scalar curvature are due to for instance \cite{GrLaw, Ros, SY, Hi74}.
Viewing the eta and rho invariant of the Dirac operator as an obstruction to the existence of 
Riemannian metrics of positive scalar curvature on compact spin manifolds was done in \cite{M92}, obstructions arising 
from covering spaces using the von Neumann trace in \cite{M98, BH} and on foliations in \cite{BH, BP}, are some amongst many 
papers on the subject.
As a corollary to Theorem \ref{thm:weit} one has,

\begin{corollary}\label{cor:pos}
In the notation of Proposition \ref{prop:getzler}, if the scalar curvature $R^g$ of the 
odd-dimensional, compact Riemannian spin manifold $Y$ is positive, then
 there exists $u_0>0$ such that for all $u \in [0, u_0]$,
the twisted Dirac operator $\dirac^\E_{uH}$ has trivial nullspace, where $H$ is a closed degree 
3 form on $Y$. In particular, 
the spectral flow $ \,{\rm sf}(\dirac^\E, \dirac^\E_{u_0 H}) $ of the family 
of  twisted Dirac operators $\{\dirac^\E_{uH} \big|\, u \in [0, u_0]\}$ 
is trivial and we deduce that,
$$
\eta(\dirac^\E_{u H}) - \eta(\dirac^\E) = - u\frac{\Rank(\E)}{(-4\pi i)^{\ell+1}} \int_Y H \wedge \widehat A(\Omega_Y^{-uH}) 
$$
Therefore for all $u \in [0, u_0]$, one has
$$
\rho_{spin}(Y, \E, uH,[g]) = \rho_{spin}(Y, \E,[g]).
$$
\end{corollary}
\medskip
Corollary \ref{cor:pos} above establishes conformal obstructions to the existence of metrics of positive scalar curvature.
By Proposition \ref{prop:sig1} and Theorem \ref{thm:weit}, we have

\begin{corollary}
Let $X$ be an even dimensional Riemannian spin manifold with spin boundary $\partial X=Y$.
Let $\widetilde H$ be a closed degree 3 form on $X$ which restricts 
to the closed degree 3-form $H$ on the boundary $Y$. 
Let $\Omega_X^{-u\widetilde H}$ denote the curvature of the Riemannian connection on $X$ with torsion tensor 
determined by $\widetilde H$
 and suppose that the scalar curvature $R^g>0$. Suppose also that the flat bundle $\widetilde
 \E$ on $X$ restricts to the flat bundle $\E$ on $Y$.
 Then 
 there exists $u_0>0$ such that for all $u \in [0, u_0]$,
$$
\Rank(\E)\, \int_X \widehat A (\Omega_X^{-u\widetilde H}) = \eta(\dirac^\E_{u H}),
$$
so that
$$
\eta(\dirac^\E_{u H}) - \eta(\dirac^\E) = \Rank(\E)\, \int_X\left[ \widehat A (\Omega_X^{-u\widetilde H}) 
- \widehat A (\Omega_X)  \right].
$$
Therefore for all $u \in [0, u_0]$, one has
$$
\rho_{spin}(Y, \E, uH,[g]) = \rho_{spin}(Y, \E,[g]).
$$

\end{corollary}

The proof of the following uses the Gromov-Lawson-Rosenberg conjecture \cite{GrLaw,Ros} and Corollary \ref{cor:pos}.
%Proposition \ref{prop:getzler},

\begin{corollary}
Let $Y$ be a closed odd dimensional Riemannian spin manifold with positive scalar curvature. 
Suppose that the Gromov-Lawson-Rosenberg conjecture  \cite{GrLaw,Ros} holds for the fundamental group $\Gamma$
of $Y$ and that $H$ is a closed 3-form on $Y$ such that 
$[H] = f^*(c)$ where $c \in H^{3}(B\Gamma)$ and $f: Y \to B\Gamma$ is a continuous map.
Then  there exists $u_0>0$ such that for all $u \in [0, u_0]$, 
$$
\eta(\dirac^\E_{u H}) = \eta(\dirac^\E).
$$
\end{corollary}

\end{document}

\title[IPositive scalar curvature and twisted Dirac complexes]
{Conformal invariants of twisted Dirac operators\\ and positive scalar curvature - Corrigendum}
%\footnote{\today}
\author{Moulay Tahar Benameur}
\address{Laboratoire et D\'epartement de Math\'ematiques,
 UMR 7122,
 Universit\'e de Metz et CNRS,
 B\^at. A, Ile du Saulcy,
 F-57045 Metz Cedex 1, 
 France}
\email{benameur@math.univ-metz.fr}

\author{Varghese Mathai}
\address{Department of Mathematics, University of Adelaide,
Adelaide 5005, Australia}
\email{mathai.varghese@adelaide.edu.au}

%
%\begin{abstract}
%We explain how to fix  an error in the statement of Theorem 2.6 in \cite{BM-JGP}. For simplicity, we have gathered the explicit details of the proof of the new statement are given in the companion paper \cite{BH-SpectralSections}. 
%\\
%\end{abstract}



\date{}
\maketitle

Thomas Schick pointed out to us that Formula (4) in \cite{BM-JGP} does not extend to the twisted case as claimed. Therefore, the sentence in \cite{BM-JGP} (right after  formula (4)):\\

``{\em{It follows that the twisted Dirac operator  $\dirac_H^{X, \maE}$ is also conformally covariant with the same weights}}''\\

is incorrect. \\

We now give the correct statement. We keep the notations of \cite{BM-JGP} and set for instance $H=\sum_j i^{j+1} H_{2j+1}$. Then with
$$
H_u := \sum_j e^{-(2j+2)u} i^{j+1} H_{2j+1} \quad \text{ and }\quad \dirac^{u g}_{H_u} = \dirac^{u g} + \hat c (H_u).
$$
and with  ${\hat c}$ denoting the Clifford representation associated with the new metric $ug$, we have  for any spinor $\psi$,  
$$
\dirac^{u g}_{H_u} (e^{-\frac{n-1}{2} u} \,\hat \psi) = e^{-\frac{n+1}{2} u} \,\widehat{\dirac^g_H\psi}.
$$
This is the content of Proposition 1.1 in \cite{BM-SpectralSections}. The conformal invariance of the twisted rho can be stated as follows (this is Theorem 1.3 in \cite{BM-SpectralSections} and should now replace  Theorem  2.6 in \cite{BM-JGP}):

\begin{theorem}[Conformal invariance of the spin rho invariant]\label{thm:indept}
The spin rho invariant $\rho_{spin}(Y,\E,H, g)$ of the twisted Dirac operator
depends only on the conformal class of the pair $(H,g)$. Said differently,
$$
\rho_{spin}(Y,\E,H_u,  u g) = \rho_{spin}(Y,\E,H, g)
$$
% and on the conformal class
%of the Hermitian metric on $\E$. 
%It is a diffeomorphism invariant of $X, \E$ and $H$.
\end{theorem}

We notice that the conformal invariance theorem was only stated in \cite{BM-JGP} and not used at all in the sequel of that paper. We also point out that the above theorem \ref{thm:indept}, although proved for general $u$, is only used in the case of constant $u$  in \cite{BM-SpectralSections}.